\documentclass[runningheads,a4paper]{llncs}

\usepackage{amssymb}
\usepackage{graphicx}

\usepackage{url}
\urldef{\mailsa}\path|{roberto.cavoretto, alessandra.derossi, hanli.qiao}@unito.it|
\urldef{\mailsb}\path|Bernhard.Quatember@uibk.ac.at|
\urldef{\mailsc}\path|Wolfgang.Recheis@i-med.ac.at|
\urldef{\mailsd}\path|Martin.Mayr@fhwn.ac.at|    
\newcommand{\keywords}[1]{\par\addvspace\baselineskip
\noindent\keywordname\enspace\ignorespaces#1}

\newcommand{\RR}{\mathbb{R}}

\begin{document}

\mainmatter  

\title{Computing Topology Preservation of RBF Transformations for Landmark-Based Image Registration}

\titlerunning{Topology Preservation in Image Registration}

%
%
\author{R. Cavoretto$^{(1)}$
\and A. De Rossi$^{(1)}$\and H. Qiao$^{(1)}$\and B. Quatember$^{(2)}$\and W. Recheis$^{(2)}$\and M. Mayr$^{(3)}$ }
\authorrunning{R. Cavoretto\and A. De Rossi\and H. Qiao\and B. Quatember\and W. Recheis\and M. Mayr }

\institute{(1) Department of Mathematics \lq\lq G. Peano\rq\rq, University of Torino,\\
Via Carlo Alberto 10, 10123 Torino, Italy\\
(2) Innsbruck Medical University,\\
Anichstrasse 35, 6020 Innsbruck, Austria\\
(3) University of Applied Science, J. Gutenberg\\
Strasse 3, 2700 Wiener Neustadt, Austria\\
\mailsa\\
\mailsb\\
\mailsc\\
\mailsd\\
\url{}}

%
%

\toctitle{Lecture Notes in Computer Science}
\tocauthor{R. Cavoretto, A. De Rossi, H. Qiao, Q. Quatember, W. Recheis, M. Mayr}
\maketitle

\begin{abstract}
In image registration, a proper transformation should be topology preserving. Especially for landmark-based image registration, if the displacement of one landmark is larger enough than those of neighbourhood landmarks, topology violation will be occurred. This paper aim to analyse the topology preservation of some Radial Basis Functions (RBFs) which are used to model deformations in image registration. Mat\'{e}rn functions are quite common in the statistic literature (see, e.g. \cite{Matern86,Stein99}). In this paper, we use them to solve the landmark-based image registration problem. We present the topology preservation properties of RBFs in one landmark and four landmarks model respectively. Numerical results of three kinds of Mat\'{e}rn transformations are compared with results of Gaussian, Wendland's, and Wu's functions.
\end{abstract}

\section{Introduction}

Over the last years, one of the largest areas of research in medical image processing has been the development of methods in image registration. The scope is to provide a point by point correspondence between the data
sets. This means to find a suitable transformation between two images, called \emph{source} and \emph{target} images, taken either at different times or from different sensors or viewpoints. The scope is to determine a transformation such that the transformed version of the source image is similar to the target one. There is a large number of applications demanding image registration, for an overview see e.g. \cite{b07Rohr01,b07Scherzer06,Zitova03}. Paper \cite{Zitova03} points out that Radial Basis Functions (RBFs) are powerful tools that could be applied to registration problem of local and global deformations. They have special good property for which values of these functions only depend on the distance of points from the center. We apply them to landmark-based registration, in which the basic idea is to determine the transformation mapping the source image onto the target image using corresponding landmarks. 

RBFs can be divided into two categories that are \emph{globally supported} and \emph{compactly supported},  respectively. In general, using \emph{globally supported} RBFs, such as thin plate spline (TPS), a single landmark pair change may influence the whole registration result but, mostly, they can keep the bending energy small. Otherwise, the Compactly Supported RBFs (CSRBFs), such as Wendland's and Wu's transformations, can circumvent this disadvantage (see \cite{Allasia14,Cavoretto12}), but usually they cannot guarantee that the bending energy is small. Papers \cite{Cavoretto13} and \cite{Derossi13} analyse different CSRBF properties for image registration which are based on global deformations.

No matter which kind of RBFs we choose to transform images, the deformed images should preserve topology. In this paper we consider the Mat\'{e}rn functions which are positive definite, and we compare the characters of topology preservation with Gaussian, Wendland's and Wu's functions (see \cite{Arad94,Fasshauer07,Gneiting10,Wendland05,Wu95}). The latter two kinds of functions are CSRBFs, whereas the former two, Mat\'{e}rn and Gaussian, are \emph{globally supported}. However, we point out that they have similar behavior and function values of them approach zeros with growing distance from their center; therefore they could be truncated as CSRBFs, being able to deal with local deformations well and allowing deformation fields to be controlled and adjusted locally using a number of landmarks points. Support size is an important index to evaluate the topology preservation property of different CSRBFs, since it modifies the influence of landmarks. In general, with small support CSRBFs can be used to deal with local image warping and with large support they can be used to deform large regions or entire images. In case of one landmark registration, we analyse the topology preservation of CSRBFs with small support size. Meanwhile in another case, a four landmarks model, we evaluate their topology properties with large support.

We arrange this paper as follows. Section 2 introduces the landmark-based registration problem for RBFs. In Section 3, three kinds of Mat\'{e}rn functions and their transformations are introduced. Section 4 compares results of topology preservation using Mat\'{e}rn transformations with results of Gaussian, Wendland's and Wu's transformation proposed by \cite{Yang11} in one landmark model. Four landmarks-based registration using transformations mentioned before are analysed in Section 5. We conclude reviewing the main results of this paper and future work in Section 6.

\section{Landmark-based Image Registration and Some RBFs}
We can interpolate the displacements defined at the landmarks using RBFs to model the deformation between a pair of objects in landmark-based image registration. To do this, we define a pair of landmark sets $\mathcal{S_N}$=\{$\textbf{x}_j \in \RR^2,j=1,2,...,N$\} and $\mathcal{T_N}$=\{$\textbf{t}_j \in \RR^2,j=1,2,...,N$\} corresponding to the \emph{source} and the \emph{target} images, respectively. The displacement can be displayed by $F_{k}: \RR^2\rightarrow \RR, k=1,2$, which has the following form

\begin{equation}\label{D}
F_k(\mathbf{x}) = \sum_{j=1}^N\alpha_{jk} \Psi\left(\parallel \mathbf{x}-\mathbf{x}_j \parallel\right),
\end{equation}
where $\Psi$ stands for a radial basis function, $\parallel \mathbf{x}-\mathbf{x}_j \parallel$ is the Euclidean distance between $\mathbf{x}$ and $\mathbf{x}_j$, and the coefficient $\alpha_{jk}$ can be calculated by two linear systems. The deformation $\textbf{f}:\RR^2\rightarrow \RR^2$ can be written as
\begin{center}
$\textbf{f}(\mathbf{x}) = \mathbf{x}+F_k(\textbf{x})$.
\end{center}

In this paper, we mainly refer to topology preservation property of RBFs that are listed here:
\begin{center}
Gaussian function: $e^{-\parallel r \parallel^2/\sigma^2}$, $\sigma>0$,
\end{center}
\begin{center}
Wendland's function $\varphi_{3,1}$: $(1-\frac{\parallel r \parallel}{c})_{+}^4(4\frac{\parallel r \parallel}{c}+1)$, $\frac{\parallel r \parallel}{c} \leqslant 1$,
\end{center}
\begin{center}
Wu's function $\psi_{1,2}$: $(1-\frac{\parallel r \parallel}{c})_{+}^4(1+4\frac{\parallel r \parallel}{c}+3(\frac{\parallel r \parallel}{c})^2+\frac{3}{4}(\frac{\parallel r \parallel}{c})^3)$, $\frac{\parallel r \parallel}{c} \leqslant 1$. 
\end{center}

Here $c$ is the support size of Wendland's and Wu's functions and $\sigma$ is the locality parameter of Gaussian. 

\section{Mat\'{e}rn Transformations}

Mat\'{e}rn functions are strictly positive define and quite common in the statistics literature \cite{Fasshauer07}. Mat\'{e}rn family has recently received a great deal of attention and has the following form \cite{Gneiting10}

\begin{equation}\label{M}
M(r \mid v,c)=\frac{2^{1-v}}{\Gamma(v)} \left(\frac{\parallel r \parallel}{c}\right)^v K_v \left(\frac{\parallel r \parallel}{c} \right).
\end{equation}

Here $K_v$ is the \emph{Modified Bessel Function of the second kind of order $v$}, $v=\beta-d/2$ and $c$ is the coefficient to determine the width or the support of functions. The Fourier transform  of the Mat\'{e}rn functions is given by the \emph{Bessel kernels}
\begin{eqnarray}
\hat{M}(w) = \left(1+ {\parallel w \parallel}^2 \right)^{-\beta} >0.
\end{eqnarray}

Therefore the Mat\'{e}rn functions are strictly positive definite, which is an important condition to ensure interpolation problem (\ref{D}) has a unique solution, and radial on $\RR^d$ for all $d<2\beta$ (see \cite{Fasshauer07}). The three specific Mat\'{e}rn functions we consider are

\begin{equation}\label{M1}
M(r \mid \frac{1}{2},c) = \frac{2^{\frac{1}{2}}}{\Gamma(\frac{1}{2})} \left(\frac{\parallel r \parallel}{c} \right)^\frac{1}{2} K_{\frac{1}{2}} \left(\frac{\parallel r \parallel}{c} \right)
                        \doteq e^{-\parallel r \parallel /c},
\end{equation}

\begin{equation}\label{M3}
M(r \mid \frac{3}{2},c) = \frac{2^{-\frac{1}{2}}}{\Gamma(\frac{3}{2})} \left(\frac{\parallel r \parallel}{c} \right)^\frac{3}{2} K_{\frac{3}{2}} \left(\frac{\parallel r \parallel}{c} \right)\doteq \left(1+\frac{\parallel r \parallel}{c} \right)  e^{-\parallel r \parallel /c},
\end{equation}

\begin{equation}\label{M5}
M(r \mid \frac{5}{2},c) = \frac{2^{-\frac{3}{2}}}{\Gamma(\frac{5}{2})} \left(\frac{\parallel r \parallel}{c} \right)^\frac{5}{2} K_{\frac{5}{2}} \left(\frac{\parallel r \parallel}{c} \right)\doteq \left(1+\frac{\parallel r \parallel}{c} + \frac{1}{3} \frac{{\parallel r \parallel}^2}{c^2}\right)  e^{-\parallel r \parallel /c}.
\end{equation}

Here $c$ is locality parameter. We title the above formulas (\ref{M1}), (\ref{M3}), (\ref{M5}) as $M_{1/2}$, $M_{3/2}$ and $M_{5/2}$, respectively. With regard to the landmark-based image registration context we give the Mat\'{e}rn transformation as follows.

\begin{definition} \label{def1}
Given a set of source landmark points $\mathcal{S_N}$=\{$\textbf{x}_j \in \RR^2,j=1,2,$ $\ldots,N$\}, and the corresponding set of target landmark points $\mathcal{T_N}$=\{$\textbf{t}_j \in \RR^2,j=1,2,\ldots,N$\}, the Mat\'{e}rn's transformation 
$\textbf{M} :\RR^2 \rightarrow \RR^2$ is such that each its component
\begin{center}
$M_k(\mathbf{x}) :\RR^2 \rightarrow \RR$,\  $k=1,2,$ 
\end{center}
assumes the following form 
\begin{equation}\label{LI1}
M_k(\mathbf{x}) = M_k(x_1,x_2) = \sum_{j=1}^N\alpha_{jk} M_{v}\left(\parallel \mathbf{x}-\mathbf{x}_j \parallel_2\right),
\end{equation}  
with $\textbf{x}$ = $(x_1,x_2)$ and $\textbf{x}_j$ = $(x_{j1},x_{j2})$ $\in \RR^2.$
\end{definition}

Through Definition \ref{def1}, we obtain transformation function $M_k(x) :\RR^2 \rightarrow \RR$ that is calculated for each $k=1,2,$ and the coefficients $\alpha_{jk}$ are to be obtained by solving two systems of linear equations.

\section{Analysis of Topology Preservation in One-Landmark Matching}

Necessary conditions to have topology preservation are continuity of the function $\textbf{H}$ and positivity of the Jacobian determinant at each point. This is achieved to get injectivity of the map \cite{Fornefett01}.

Suppose that the source landmark $\textbf{p}$ is shifted by $\Delta_x$ along the $x$-axis direction and by $\Delta_y$ along the $y$-axis direction to the target landmark $\textbf{q}$. The coordinates of transformation are
\begin{center}
$H_1(\textbf{x})=x+\Delta_x\Phi(||\textbf{x}-\textbf{p}||)$,  
\end{center}  
\begin{center}
$H_2(\textbf{x})=y+\Delta_y\Phi(||\textbf{x}-\textbf{p}||)$,
\end{center}
where $\Phi$ is any RBF.

Requiring the determinant of the Jacobian is positive, we obtain
\begin{equation}
\det(J(x,y))=1+\Delta_x\frac{\partial \Phi}{\partial x}+\Delta_y\frac{\partial \Phi}{\partial y}>0,
\end{equation}
i.e.
\begin{center}
$\Delta_x\frac{\partial \Phi}{\partial x}+\Delta_y\frac{\partial \Phi}{\partial y}>-1$,
\end{center}
or, equivalently,
\begin{center}
$\Delta_x\frac{\partial \Phi}{\partial r}\cos\theta+\Delta_y\frac{\partial \Phi}{\partial r}\sin\theta>-1$,
\end{center}
where $\Phi$ stands for $\Phi(||\textbf{x}-\textbf{p}||)$ and $r=||\textbf{x}-\textbf{p}||$. 

If we set $\Delta=\mbox{max}(\Delta_x,\Delta_y)$, the value of $\theta$ minimizing the determinant in 2D is $\frac{\pi}{4}$; thus we get
\begin{equation}\label{WC1}
\Delta\frac{\partial\Phi}{\partial r}>-\frac{1}{\sqrt{2}}.
\end{equation}
With the condition (\ref{WC1}) one can show that all principal minors of the Jacobian are positive. It follows that the transformations defined by equation (\ref{LI1}) preserve the topology if (\ref{WC1}) holds. The minimum of $\frac{\partial\Phi}{\partial r}$ depends on the localization parameter and therefore on the support size of the parameter $c$ of Mat\'{e}rn functions.

In the next subsections we compute the minimum support size of locality parameter, for which (\ref{WC1}) is satisfied, of three functions of Mat\'{e}rn family. 

\subsection{Mat\'{e}rn $M_{1/2}$} 

Now we are considering the Mat\'{e}rn function (\ref{M1}). Clearly, we cannot obtain the minimum value of $c$ as the above mentioned method, so we search for it through numerical experiments. Here we have
\begin{equation}\label{$M_1.1$}
\frac{\partial\Phi}{\partial r}=-\frac{1}{c} e^{-r/c},
\end{equation}
while the value of $r$ minimizing (\ref{$M_1.1$}) is given by $r=0.25c$. Then computing (\ref{$M_1.1$}) and substituting its result in (\ref{WC1}), we obtain
\begin{equation}
c>\frac{\sqrt{2}}{e^{1/4}}\Delta\approx 1.1\Delta.
\end{equation}

\subsection{Mat\'{e}rn $M_{3/2}$} 
According to Mat\'{e}rn function (\ref{M3}), we find the minimum value of $c$ satisfying (\ref{WC1}). Now we get
\begin{equation}\label{$M_3.1$}
\frac{\partial\Phi}{\partial r}=-\frac{r}{c^2}e^{-r/c},
\end{equation}
while the value of $r$ minimizing (\ref{$M_3.1$}) is given by $r = c$. Evaluating (\ref{$M_3.1$}) at $r = c$  and substituting its result in (\ref{WC1}), we obtain
\begin{equation}
c>\frac{\sqrt{2}\Delta}{e}\approx 0.52\Delta.
\end{equation}

\subsection{Mat\'{e}rn $M_{5/2}$} 
Considering Mat\'{e}rn function (\ref{M5}), similarly as Mat\'{e}rn function (\ref{M3}), we have
\begin{equation}\label{$M_5.1$}
\frac{\partial\Phi}{\partial r}=-\bigg(\frac{r}{3c^2}+\frac{r^2}{3c^2}\bigg)e^{-r/c},
\end{equation}
while the value of $r$ minimizing (\ref{$M_5.1$}) is given by $r = \frac{\sqrt{5}+1}{2}c$. Calculating (\ref{$M_5.1$}) and substituting its result in (\ref{WC1}), we obtain
\begin{equation}
c>\frac{(2\sqrt{2}+\sqrt{10})\Delta}{3e^{(\sqrt{5}+1)/2}}\approx 0.3960\Delta.
\end{equation}

\subsection{Analysis of the Results}
Table \ref{tab1} summarizes the minimum support sizes of locality parameter for $M_{1/2}$, $M_{3/2}$ and $M_{5/2}$, which are compared with Gaussian, Wendland's and Wu's functions (see \cite{Fornefett01,Yang11}). The advantage of having small supports is that the influence area of each landmark turns out to be small. This allows us to have a greater local control. From Table \ref{tab1}, we can see that Mat\'{e}rn functions have smaller supports, especially $M_{5/2}$ function. This means that in one landmark model, the deformed field of $M_{5/2}$ is the smallest among these six transformations. 

\begin{table}
\caption{Minimum support size for various RBFs, where $c=2\sigma$ and $d=2$.}
\begin{center}\label{tab1}
\begin{tabular}{cccccc}
\hline
$Gaussian$ & $\varphi_{3,1}$ & $\psi_{1,2}$ & $M_{1/2}$ & $M_{3/2}$ & $M_{5/2}$ \\
\hline
$\sigma>1.21\Delta$ & $c>2.98\Delta$ & $c>2.80\Delta$ & $c>1.10\Delta$ & $c>0.52\Delta$ & $c>0.3960\Delta$\\ 
\hline
\end{tabular}
\end{center}
\end{table}

\subsection{Numerical Results}

\begin{figure}
\begin{minipage}{62mm}
\includegraphics[width=6.2cm]{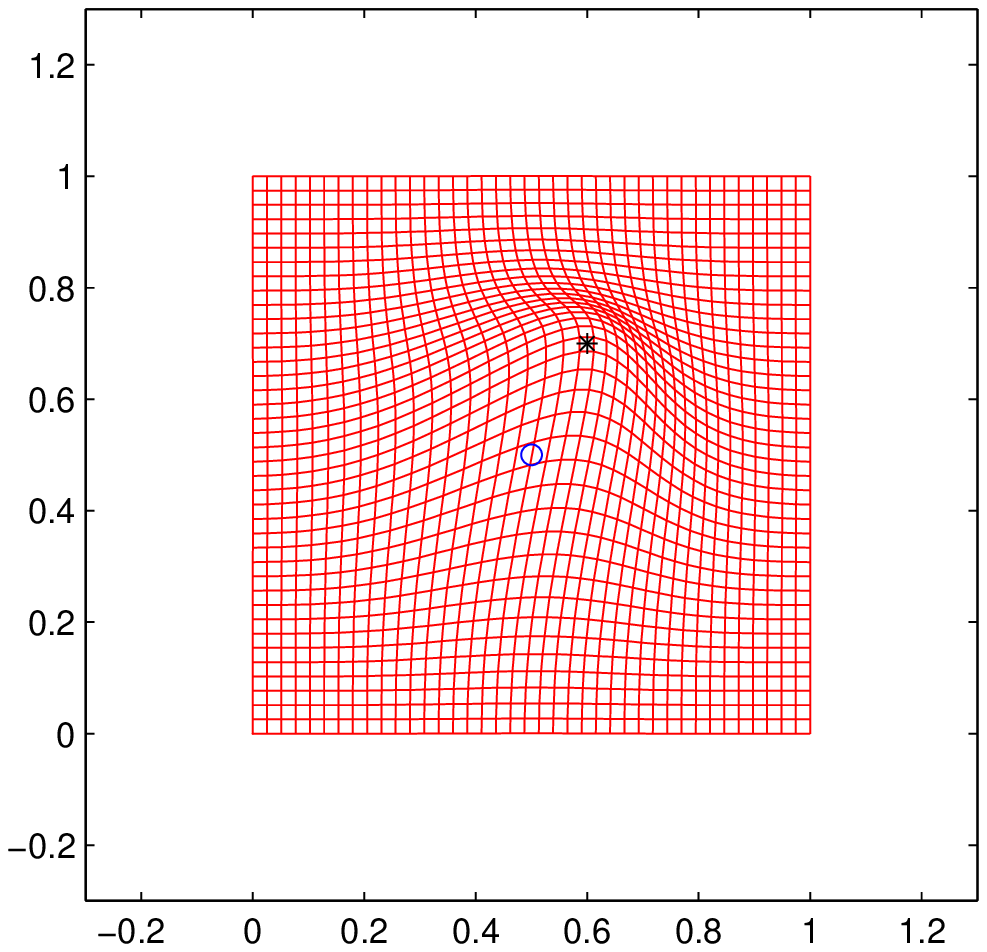}
\centerline{(a) Wendland $\varphi_{3,1}$, $c = 0.6$}
\end{minipage}
\begin{minipage}{62mm}
\includegraphics[width=6.2cm]{{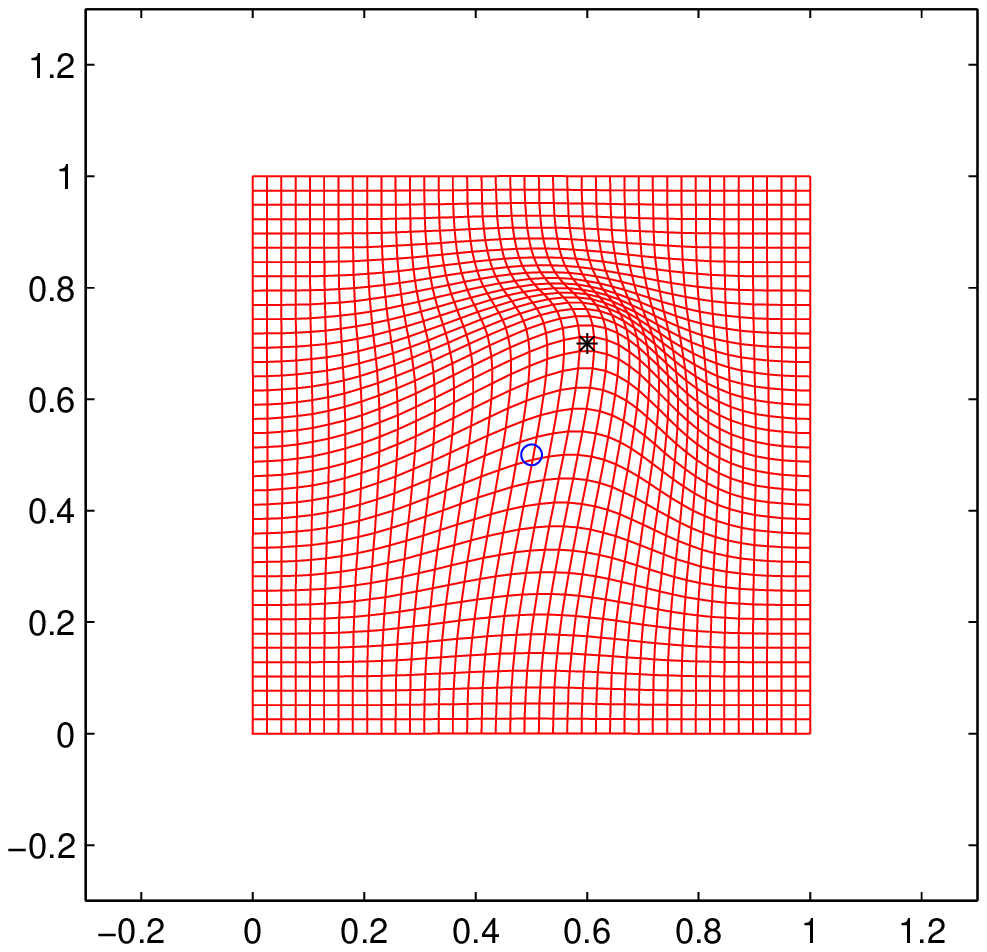}}
\centerline{(b) Wu $\psi_{1,2}$, $c = 0.58$}
\end{minipage}
\begin{minipage}{62mm}
\includegraphics[width=6.2cm]{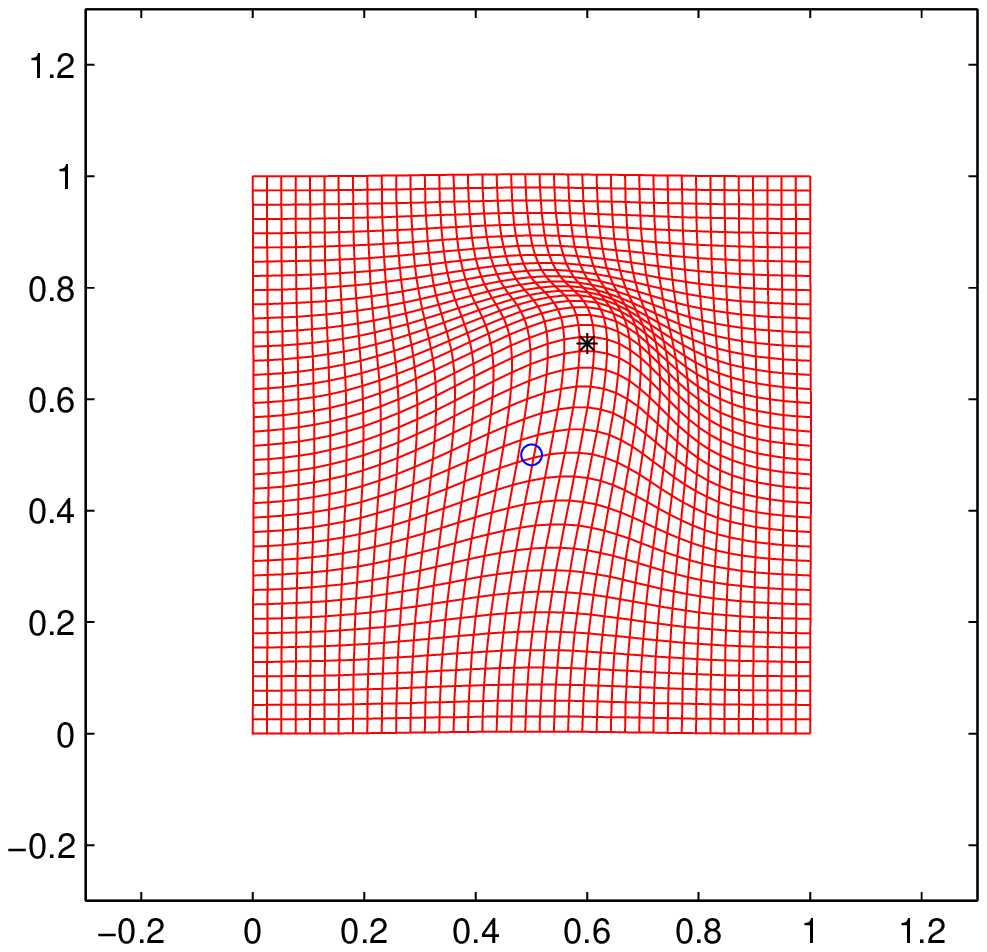}
\centerline{(c) Gaussian $\sigma = 0.25$}
\end{minipage}
\begin{minipage}{62mm}
\includegraphics[width=6.2cm]{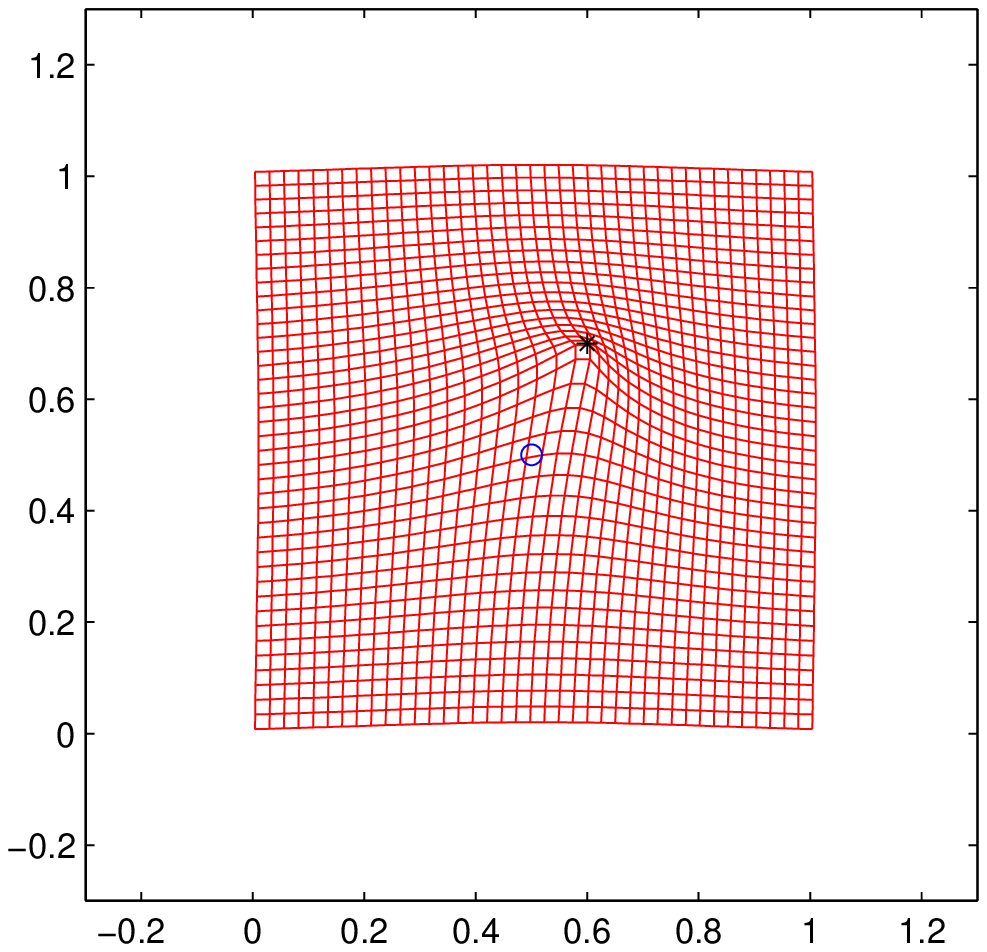}
\centerline{(d) Mat\'{e}rn, $M_{1/2}$ $c = 0.22$}
\end{minipage}
\begin{minipage}{62mm}
\includegraphics[width=6.2cm]{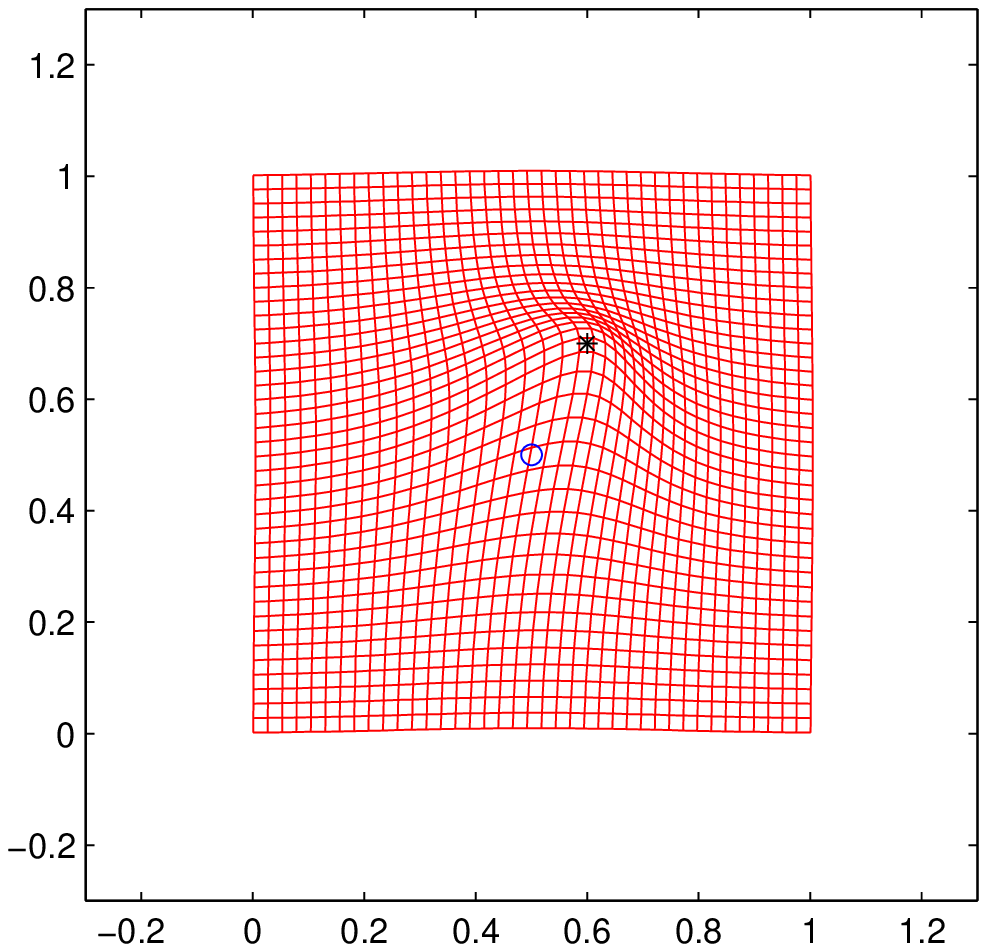} 
\centerline{(e) Mat\'{e}rn, $M_{3/2}$ $c = 0.105$}
\end{minipage}
\begin{minipage}{62mm}
\includegraphics[width=6.2cm]{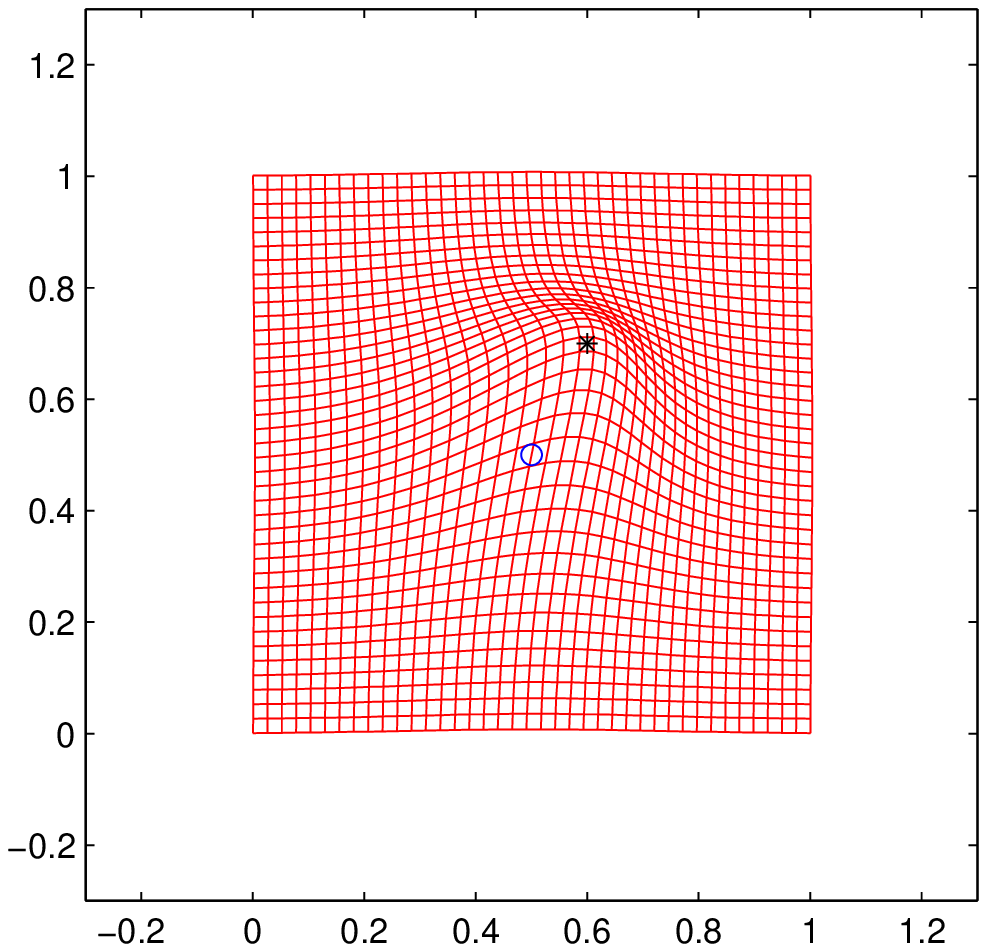}
\centerline{(f) Mat\'{e}rn, $M_{5/2}$ $c = 0.08$}
\end{minipage}
\caption{Deformation results of one-landmark matching using minimum locality parameters satisfying the topology preservation condition. The source landmark is marked by a circle ($\circ$), while the target one by a star ($\ast$).}
\label{fig_1}
\end{figure}

In this section, we report the numerical experiments obtained on a grid $[0,1]\times[0,1]$ and compare then the distortion outcomes of the grid in the shift case of the landmark $\{(0.5,0.5)\}$ in $\{(0.6,0.7)\}$. In Figure \ref{fig_1} we show results assuming as a support size the minimum $c$ and $\sigma$ such that (\ref{WC1}) is satisfied, with $\Delta=0.2$.

Figure \ref{fig_1} shows that, for the minimum value of $c$ and $\sigma$, all transformations can preserve topology well. In this case, $M_{5/2}$ transformation has the smallest deformed field around the landmark while Wendland's, Wu's and Gaussian has relatively larger field, as outlined in Table \ref{tab1}. Nevertheless, we can see that the whole images are slightly deformed. In other words, if the topology preservation condition (\ref{WC1}) is not satisfied, the transformed image is deeply misrepresented above all around the shifted point.

\section{Topology Preservation for More Extended Deformations}

If we consider much large supports which are able to cover whole domain, the influence of each landmark extends on the entire image, thus generating global deformations. In the following we compare topology preservation properties for globally supported transformations. For this aim, we consider four inner landmarks in a grid, located so as to form a rhombus at the center of the figure, and we suppose that only the lower vertex is downward shifted of $\Delta$ \cite{Yang11}. The landmarks of source and target images are
$P=\{(0,1),(-1,0),(0,-1),(1,0)\}$ and $Q=\{(0,1),(-1,0),(0,-1-\Delta),(1,0)\}$, respectively, with $\Delta>0$.

Let us now consider components of a generic transformation $\textbf{H}:\RR^2\rightarrow \RR^2$ obtained by a transformation of four points $P_1$, $P_2$, $P_3$ and $P_4$, namely
\begin{equation}
H_1(\textbf{x})=x+\sum_{i=1}^4 c_{1,i} \Phi(||\textbf{x}-P_i||),
\end{equation}
\begin{equation}
H_2(\textbf{x})=y+\sum_{i=1}^4 c_{2,i}\Phi(||\textbf{x}-P_i||).
\end{equation}
The coefficients $c_{1,i}$ and $c_{2,i}$ are obtained so that the transformation sends $P_i$ to $Q_i$, with $i=1,\dots ,4$. To do that, we need to solve two systems of four equations in four unknowns, whose solutions are
\begin{center}
$c_{1,1}=0,\ \ \ c_{1,2}=0,\ \ \ c_{1,3}=0,\ \ \ c_{1,4}=0$,
\end{center}
and
\begin{equation}\label{condition4_1}
c_{2,1} = \frac{\beta^2+\beta-2\alpha^2}{(1-\beta)[(1+\beta)^2-4\alpha^2]}\Delta, \ \ \ c_{2,2} = \frac{\alpha}{(1+\beta)^2-4\alpha^2}\Delta,
\end{equation}
\begin{equation}\label{condition4_2}
c_{2,3} = -\frac{1+\beta-2\alpha^2}{(1-\beta)[(1+\beta)^2-4\alpha^2]}\Delta, \ \ \ c_{2,4} = c_{2,2},
\end{equation}
where $\alpha=\Phi\left(\frac{\sqrt{2}}{c}\right)$ and $\beta=\Phi\left(\frac{2}{c}\right).$ For simplicity, we denote $\Phi_1=\Phi(||(x,y)-P_1||/c)$, $\Phi_2=\Phi(||(x,y)-P_2||/c)$, $\Phi_3=\Phi(||(x,y)-P_3||/c)$ and $\Phi_4=\Phi(||(x,y)-P_4||/c)$.

The determinant of the Jacobian matrix is
\begin{equation}\label{condition4_3}
\det\left(J(x,y)\right)=1+\sum_{i=1}^4 c_{2,i}\frac{\partial\Phi_i}{\partial y}.
\end{equation}
The minimum Jacobian determinant is obtained at position $(0,y)$, with $y>1$. In the following, we analyse the value of the Jacobian determinant at $(0,y)$, with $y>1$ for different RBFs. Since the support $c$ is very large, in order to have a global transformation, we consider $||\cdot||/c$ to be infinitesimal and omit terms of higher order. 

\subsection{Mat\'{e}rn $M_{1/2}$} \label{sec:4.3}
We approximate Mat\'{e}rn function $M_{1/2}$ as follows
\begin{equation}
M_{1/2}(r)=e^{-r}\approx 1-r+\frac{r^2}{2},
\end{equation}
while its first derivative is $M_{1/2}'(r)\approx -1+r$. Now, we approximate $\alpha$ and $\beta$ as $\alpha=\Phi\left(\frac{\sqrt{2}}{c}\right)\approx1-\frac{\sqrt{2}}{c}+\frac{1}{c^2}$, $\beta=\Phi\left(\frac{2}{c}\right)\approx1-\frac{2}{c}+\frac{2}{c^2}$. Based on the approximated $\alpha$ and $\beta$, and according to (\ref{condition4_1})--(\ref{condition4_3}) we obtain the Jacobian determinant of $M_{1/2}$ transformation as the following form: 
\begin{equation}
\det\left(J(0,y)\right) \approx 1-2.4142\Delta\left(-1+\frac{y}{\sqrt{y^2+1}}\right).
\end{equation}

\subsection{Mat\'{e}rn $M_{3/2}$} \label{sec:4.4}
The approximation of Mat\'{e}rn function $M_{3/2}$ is
\begin{equation}
M_{3/2}(r)=(1+r)e^{-r}\approx 1-\frac{r^2}{2}+\frac{r^3}{2},
\end{equation}
its first derivative is $M_{3/2}'(r)=-re^{-r}\approx -r+r^2$, $\alpha$ and $\beta$ can be approximated by $\alpha=\Phi\left(\frac{\sqrt{2}}{c}\right)\approx 1-\frac{1}{c^2}+\frac{\sqrt{2}}{c^3}$, $\beta=\Phi\left(\frac{2}{c}\right)\approx 1-\frac{2}{c^2}+\frac{4}{c^3}$. Then, the $M_{3/2}$ Jacobian determinant can be formed:
\begin{equation}
\det\left(J(0,y)\right) \approx 1-1.7071\Delta\left(y^2+1-y\sqrt{y^2+1}\right).
\end{equation}

\subsection{Mat\'{e}rn $M_{5/2}$} \label{sec:4.5}
Mat\'{e}rn function $M_{5/2}$ can be approximated as
\begin{equation}
M_{5/2}(r)=(1+r+\frac{r^2}{3})e^{-r}\approx 1-\frac{r^2}{6}+\frac{2}{9}r^3,
\end{equation}
its first derivative is $M_{5/2}'(r)=-\frac{1}{3}r+\frac{2}{3}r^2$. Approximation of $\alpha$ and $\beta$ is $\alpha=\Phi\left(\frac{\sqrt{2}}{c}\right)\approx 1-\frac{1}{3c^2}+\frac{4\sqrt{2}}{9c^3}$, $\beta=\Phi\left(\frac{2}{c}\right)\approx 1-\frac{2}{3c^2}+\frac{16}{9c^3}$. Under formulas (\ref{condition4_1})--(\ref{condition4_3}), we can get the Jacobian determinant function of $M_{5/2}$, i.e.
\begin{equation}
\det\left(J(0,y)\right) \approx 1-2.5607\Delta\left(y^2+1-y\sqrt{y^2+1}\right).
\end{equation}

\subsection{Analysis of the Results}
\begin{figure}
\centering
\includegraphics[width=7.cm]{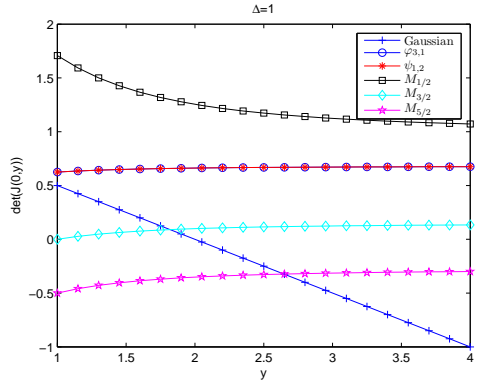}
\caption{Value of $\det(J(0,y))$, with $y>1$, by varying RBFs.}
\label{fig_2}
\end{figure}

The obtained results, if compared with the ones acquired by the work \cite{Yang11}, show same values of $\det(J(0,y))$, with $y>1$, when one uses CSRBF transformations based on Wendland's and Wu's functions. This indicates that functions $\varphi_{3,1}$, $\psi_{2,1}$ have a very similar behavior. The equations obtained in \cite{Yang11} using Wendland's and Wu's functions guarantee the Jacobian determinant positivity for any $y>1$. When one uses Gaussian, $M_{1/2}$, $M_{3/2}$ and $M_{5/2}$ functions, we can find $M_{1/2}$ and $M_{3/2}$ functions guarantee the positivity of the Jacobian determinant for any $y>1$. Using $M_{1/2}$ function, the Jacobian determinant is the closest to 1. This means that it is the best transformation in this case. While that obtained for the Gaussian, still in \cite{Yang11}, presents a negative determinant for some values of $y$, as shown in Figure \ref{fig_2}. We also see that $M_{5/2}$ function always present a negative determinant for different $y$. Therefore, $\varphi_{3,1}$, $\psi_{2,1}$, and the first two kinds of Mat\'{e}rn functions ensure more easily the topology preservation, unlike the Gaussian and $M_{5/2}$ function. We can conclude that the two functions do not lead to good results in case of higher landmarks density, i.e. when distance among landmarks is very little. Moreover, each of them influences the whole image, which might produce a topology violation.

\subsection{Numerical Results}

\begin{figure}
\begin{minipage}{62mm}
\includegraphics[width=6.2cm]{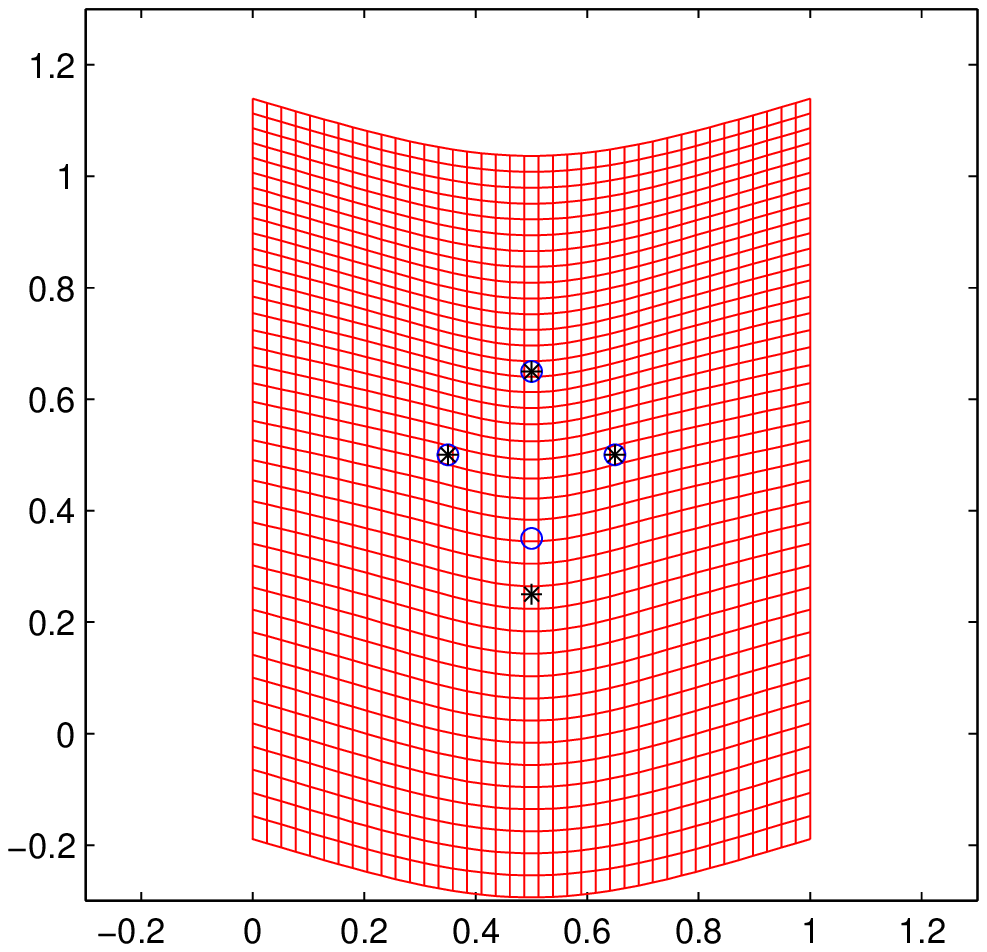}
\centerline{(a) Wendland $\varphi_{3,1}$, $c = 100$}
\end{minipage}
\begin{minipage}{62mm}
\includegraphics[width=6.2cm]{{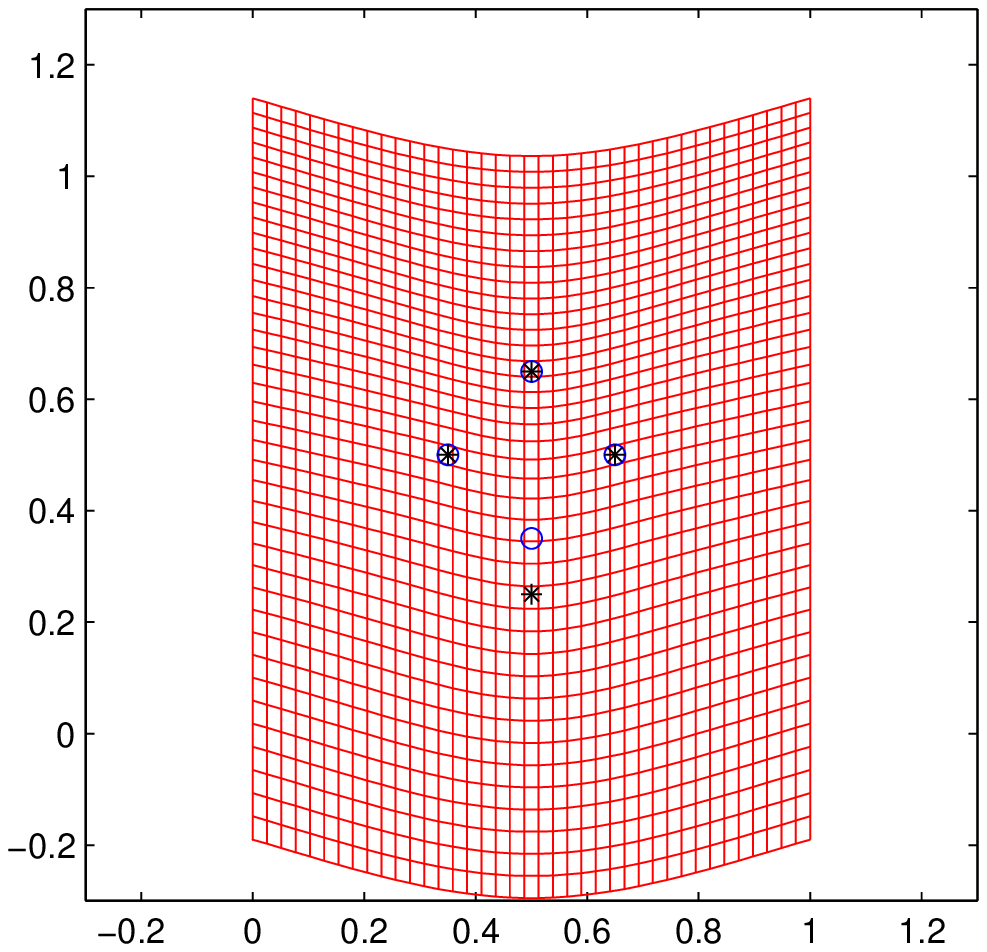}}
\centerline{(b) Wu $\psi_{1,2}$, $c = 100$}
\end{minipage}
\begin{minipage}{62mm}
\includegraphics[width=6.2cm]{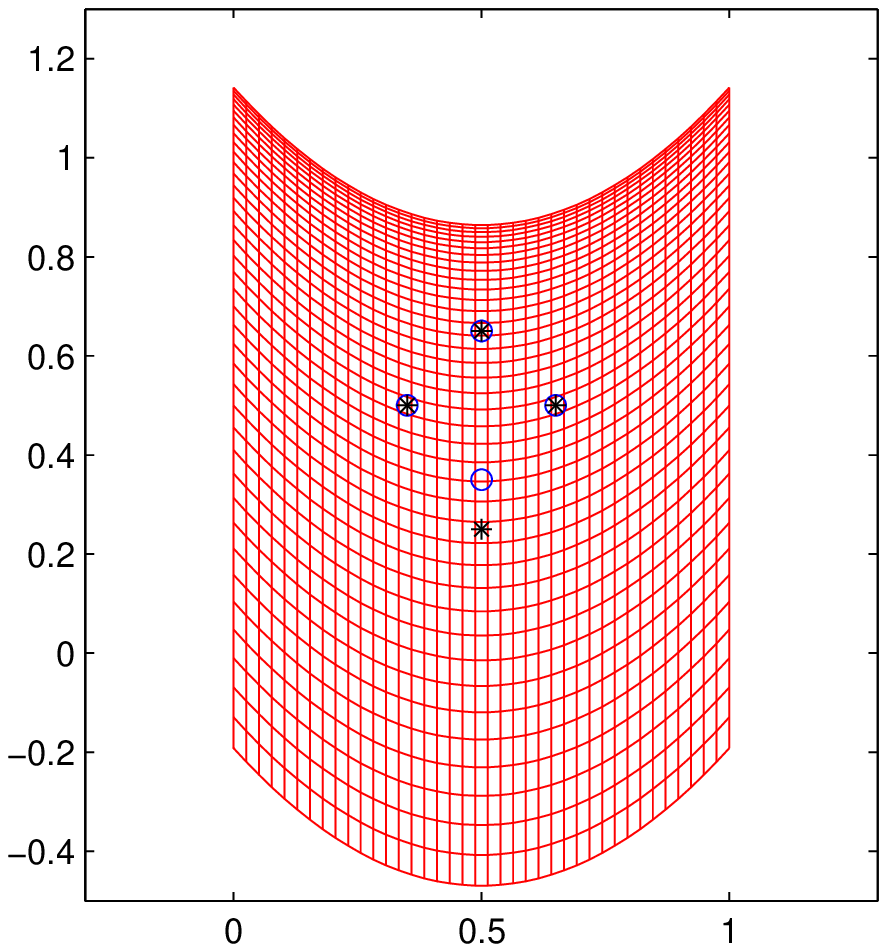}
\centerline{(c) Gaussian $\sigma = 50$}
\end{minipage}
\begin{minipage}{62mm}
\includegraphics[width=6.2cm]{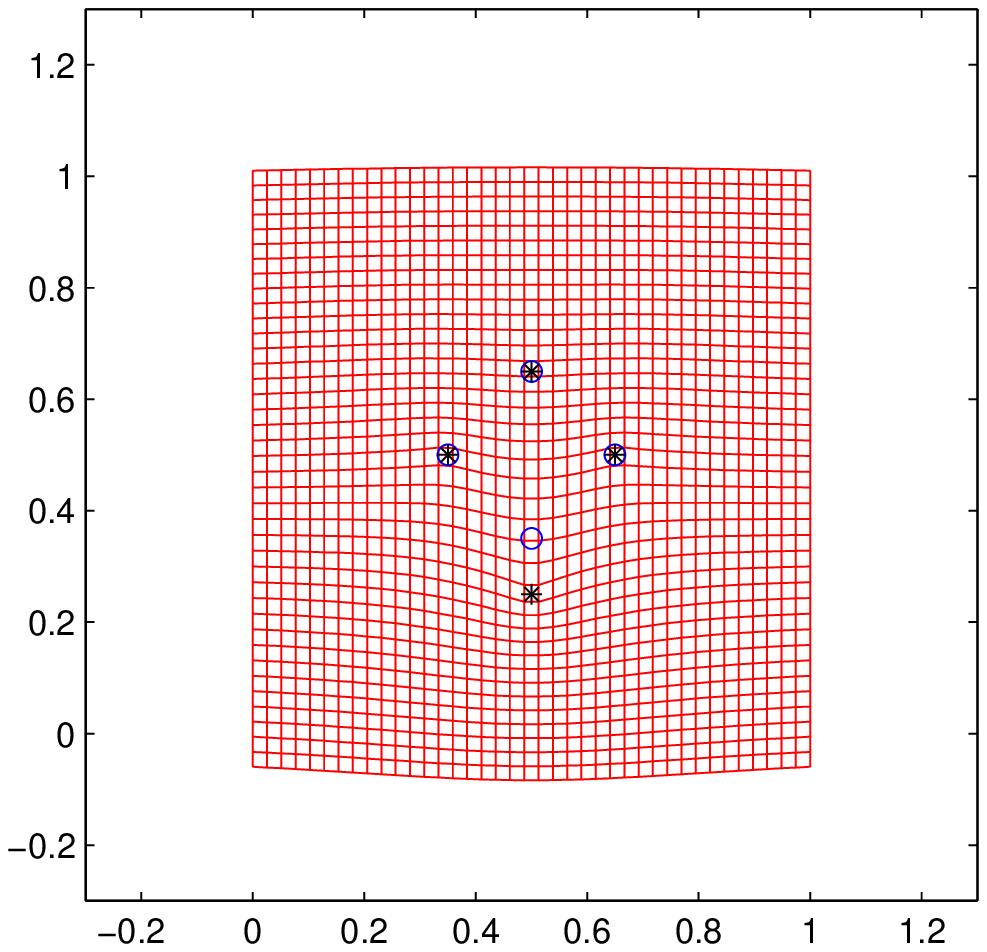}
\centerline{(d) Mat\'{e}rn, $M_{1/2}$ $c = 100$}
\end{minipage}
\begin{minipage}{62mm}
\includegraphics[width=6.2cm]{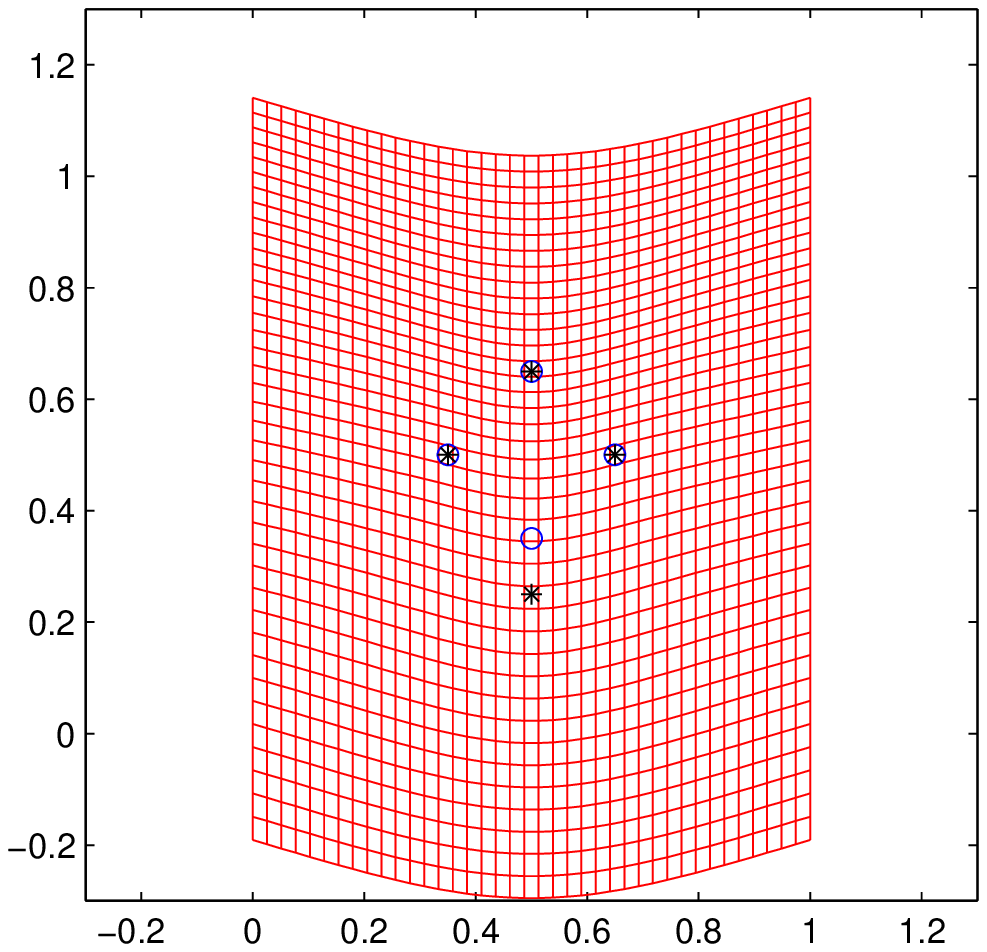} 
\centerline{(e) Mat\'{e}rn, $M_{3/2}$ $c = 100$}
\end{minipage}
\begin{minipage}{62mm}
\includegraphics[width=6.2cm]{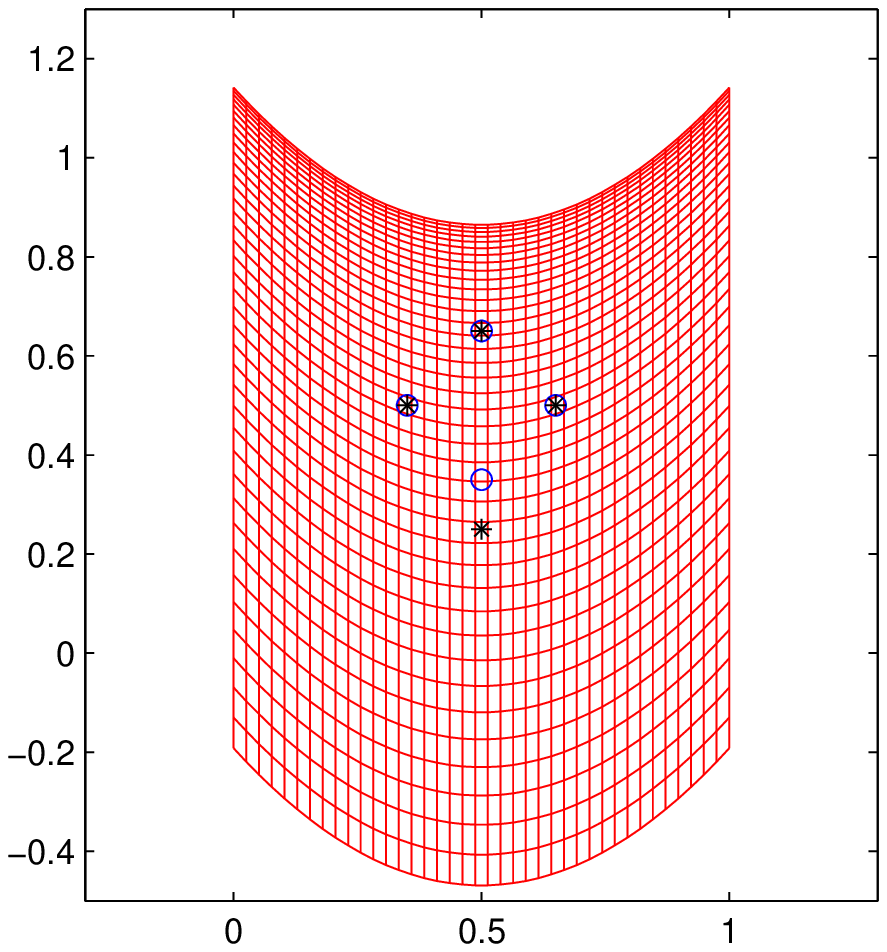}
\centerline{(f) Mat\'{e}rn, $M_{5/2}$ $c = 100$}
\end{minipage}
\caption{Deformation results of four landmarks; the source landmarks are marked by a circle ($\circ$), while the target ones by a star ($\ast$).}
\label{fig_3}
\end{figure}

Let us consider $[0,1]\times[0,1]$ and compare results obtained by its distortion, which is created by the shift of one of the four landmarks distributed in rhomboidal position. The source landmarks are $\{(0.5,0.65)$, $(0.35,0.5)$, $(0.65,0.5)$, $(0.5,0.35)\}$ and are respectively transformed in the following target landmarks $\{(0.5,0.65),$ $(0.35,0.5)$, $(0.65,0.5)$, $(0.5,0.25)\}$.
Taking $\sigma=50$ and $c=100$ as support size, we obtain Figure \ref{fig_3}. 

In agreement with theoretical results, Figure \ref{fig_3} confirms that Gaussian and $M_{5/2}$ function turns out to be those which worse preserve topology, whereas all other functions present very similar deformations. In particular, the $M_{1/2}$ function provides the best transformation. In fact, although support size is large, the deformed field at landmarks is very small in $M_{1/2}$ transformation.

\section{Conclusions and Future work}

We evaluated the topology preservation property of three kinds of Mat\'{e}rn functions with Gaussian, Wendland's and Wu's functions in one and four landmarks cases, respectively. No matter in which case, $M_{1/2}$ transformation showed us the best advantage. We must note that although in one landmark model $M_{1/2}$ deformation has not the smallest locality parameter among the six transformations, it guarantees images only deformed around the landmark instead of in the whole images which is a goodness in local deformation. In the next step we will evaluate topology property and other characters of these Mat\'{e}rn functions in real life cases, such as $x$-ray images of patients. Also, we are going to check the outcomes of Mat\'{e}rn functions in case of a large number of landmarks.

Observing the formula (\ref{M}), we found that when $v$ is large enough, Mat\'{e}rn function can be approximated as Gaussian function. Also in this paper, we can see that $M_{5/2}$ function has a character more similar to Gaussian than $M_{1/2}$ and $M_{3/2}$. Therefore another work might be analysing and comparing behaviors of Mat\'{e}rn and Gaussian functions in landmark-based image registration.

\subsubsection*{Acknowledgments.} The second author acknowledges financial support from the GNCS--INdAM. 

\end{document}